\documentclass[twocolumn]{autart}    
\usepackage [latin1]{inputenc}                                     
\usepackage{matlab-prettifier}                                     
\usepackage{amssymb}
\usepackage{enumitem}
\usepackage{graphicx}
\usepackage{CJK}
\usepackage{indentfirst}
\usepackage{amsmath}
\usepackage{color}
\usepackage{subfig}
\usepackage{wrapfig}
\newtheorem{corollary}{Corollary}
\newtheorem{assumption}{Assumption}
\newtheorem{proposition}{Proposition}

\newtheorem{lemma}{Lemma}
\newtheorem{theorem}{Theorem}

\usepackage{cases}
\begin{document}

\begin{frontmatter}

\title{Observer-Based Adaptive Scheme for Fixed-Time Frequency Estimation of Biased Sinusoidal Signals} 
\thanks[footnoteinfo]{Corresponding author: Shang Shi.}


\author[hehai,footnoteinfo]{Shang Shi}\ead{shishangshang@yahoo.com;shishangshang@foxmail.com}, \
\author[nanligong]{Huifang Min}, \ 
\author[zhenjiang]{Shihong Ding}
\address[hehai]{College of Energy and Electrical Engineering, Hohai University, Nanjing, Jiangsu 211100, PRChina}
\address[nanligong]{School of Automation, Nanjing University of Science and Technology, Nanjing, Jiangsu 210094, PRChina}
\address[zhenjiang]{
School of Electrical and Information Engineering, Jiangsu University, Zhenjiang, Jiangsu 212013, PRChina.}

\begin{abstract}
In this technical communique, we propose a novel observer-based adaptive scheme to deal with the parameter estimation problem of biased sinusoidal signals. Different from the existing  adaptive frequency estimation scheme, the proposed scheme can achieve fixed-time frequency estimation, whose convergence time is independent of the initial errors. Simulation example with different initial values shows the effectiveness of the theoretical result.
\end{abstract}

\end{frontmatter}
\section{Introduction}
The problem of frequency estimation for sinusoidal signals is a very important and  fundamental issue in both theoretical and practical applications, such as the rejection of periodic disturbance ({\color{blue}{Shi, Xu, Gu, \& Zhang, 2019}}) and the control of power systems ({\color{blue}{Rao, Soni, Sinha, \& Nasiruddin, 2019}}). Since then, a large number of approaches have been proposed to solve the frequency estimation problem including Kalman filters ({\color{blue}{Hajimolahoseini, Taban, \& Soltanian-Zadeh, 2012}}), adaptive notch filters ({\color{blue}{Hsu, Ortega, \& Damm, 1999}}), time domain-based methods ({\color{blue}{Angrisani, D'Apuzzo, Grillo, Pasquino, \& Moriello, 2014}}), algebraic identification ({\color{blue}{Trapero, Sira-Ram¨ªrez, \& Batlle, 2007}}), adaptive phase-locked-loop approaches ({\color{blue}{Karimi-Ghartemani, \& Ziarani,  2004}}) and  state-variable filtering techniques ({\color{blue}{Pyrkin, Bobtsov, Efimov, \& Zolghadri, 2011}}).

Another important class of algorithms used in the frequency estimation is the so-named adaptive observer approach. By modeling sinusoidal signals as  observable linear systems where the frequency is treated as unknown parameter, an observer with adaptive parameter can be designed to achieve frequency identification ({\color{blue}{Marino, \& Tomei, 2002; Xia, 2002; Hou, 2012}}). In {\color{blue}{Chen, Pin, Ng, Hui and Parisini (2017)}}, the adaptive observer-based estimation scheme was used to estimate the frequency of single sinusoidal signals with structured and unstructured measurement disturbances. The frequency estimation problem for more complicated multiple sinusoidal signals with bounded perturbations on the measurement was addressed in {\color{blue}{Pin, Wang, Chen and Parisini (2019)}} by using the adaptive observer approach.

Note that most of the existing results in the literatures can only achieve asymptotic frequency estimation. In the recent work {\color{blue}{Pin, Chen and Parisini (2017)}}, by using a volterra operator combined with a second-order sliding mode-based adaptation law, a new volterra operator based adaptive frequency estimator was developed, which can achieve finite-time frequency estimation of biased sinusoidal signals. In {\color{blue}{Li, Fedele, Pin and Parisini (2016)}}, the algorithm in {\color{blue}{Pin et al. (2017)}} was extended for the parameter estimation of a biased and damped sinusoidal signal. Inspired by the work of {\color{blue}{Pin et al. (2017)}}, the finite-time estimation problem of multiple biased and damped sinusoidal signals was solved in the most recent paper {\color{blue}{Chen, Li, Pin, Fedele and Parisini (2019)}}.

It can be clearly seen that the convergence time of the conventional adaptive estimator is dependent on the initial estimation errors and will grow as the initial errors grow. Although finite-time adaptive estimator has a faster convergence speed, the settling time still depends on the initial estimation errors. As an exception, the finite-time estimator in {\color{blue}{Chen et al. (2019)}} was designed by using an algebraic method. To overcome this drawback, the notation of fixed-time stability was proposed ({\color{blue}{Andrieu, Praly, \& Astolfi, 2008; Polyakov, 2012}}). Based on this notation, many remarkable results have been developed ({\color{blue}{Polyakov, 2012}}). However, to the best of the authors' knowledge, the results about fixed-time adaptive frequency estimator haven't been reported in the literature.

In this technical communique, a novel fixed-time estimation algorithm is proposed for estimating of biased sinusoidal signals. To design the estimator, an fixed-time observer-based adaptive scheme is developed. Different from the existing asymptotic and finite-time adaptive estimators, the convergence time of the proposed algorithm is bounded by a fixed time which is independent of the initial errors. This is also the main contribution of the technical communique.

\emph{Notation:} Throughout the technical communique, we define $\lfloor x\rceil^\alpha=|x|^\alpha sign(x),\forall \alpha>0, x\in \mathbb{R}$.
\section{Problem Formulation and Preliminaries}
The biased sinusoidal signal considered in this technical communique is presented as follows:
\begin{eqnarray}\label{2.1}
y(t)=A+Bsin(\phi(t)), \quad \dot{\phi}(t)=w,\quad \phi(0)=\phi_0,
\end{eqnarray}
where $y(t)$ is measurable with its derivatives $y^{(i)}(t),i\in \mathbb{N}+$ unmeasurable; $A\in\mathbb{R}_{>0},B\in\mathbb{R}_{>B_{min}},w\in\mathbb{R}_{0\cup[w_{min},+\varpropto)}$ and $\phi_0\in\mathbb{R}$ are unknown offset, amplitude, angular frequency and initial phase shift with $B_{min},w_{min}\in\mathbb{R}_{>0}$ can chosen arbitrary small.


\begin{assumption}
There exists a known positive constant $L$ and a known positive integer $m\geq4$ such that the $m$-order derivative of $y(t)$ satisfies $|y^{(m)}(t)|\leq L$.
\end{assumption}

The objective here is to estimate the parameters $w$ in a fixed-time independent of initial condition.

Firstly, an arbitrary order differentiator designed in {\color{blue}{Angulo, Moreno and Fridman (2013)}} will be used here as an observer to estimate the signal and its derivatives:
\begin{align}\label{2.2}
\dot{z}_{i}&=-\kappa_{i}\theta(t)\lfloor\tilde{z}_{1}\rceil^{\frac{m-i}{m}}-k_{i}(1-\theta(t))\lfloor\tilde{z}_{1}\rceil^{\frac{m+\alpha i}{m}}+z_{i+1},\nonumber\\
i&=1,2,\cdots,m-1,\nonumber\\
\dot{z}_{m}&=-\kappa_{m}\theta(t)sign(\tilde{z}_{1})-k_{m}(1-\theta(t))\lfloor\tilde{z}_{1}\rceil^{1+\alpha},
\end{align}
where $\tilde{z}_{1}=z_{1}-y(t)$; $\theta(t)=\frac{sign(t-T_{u})+1}{2}$ with arbitrarily chosen $T_{u}>0$; $\{\kappa_{i},k_{i}\}_{i=1}^{m}$ and $\alpha>0$ are design parameters selected the same as that in Theorem 1 of {\color{blue}{Angulo et al. (2013)}}; the states $z_{1},z_{2},\cdots,z_{m}$ are the estimations of $y(t),y^{(1)}(t),\cdots,y^{(m-1)}(t)$.

\begin{lemma}({\color{blue}{Angulo, et al., 2013}})
For the biased sinusoidal signal $y(t)$ under Assumption 1 and the observer (\ref{2.2}), all the signals of system (\ref{2.2}) are bounded and there exists a time $T_{1}$ independent of initial condition such that $z_{i}=y^{(i-1)}(t),\forall t\geq T_{1}$ is satisfied.
\end{lemma}

\begin{lemma}({\color{blue}{Polyakov, 2012}})
For system $\dot{x}=f(x)$ with $f(0)=0$, if there exists a continuous radially unbounded and positive definite function $V(x)$ such that $\dot{V}(x)\leq-\alpha V^{1+\frac{1}{\mu}}-\beta V^{1-\frac{1}{\mu}}$ with $\alpha,\beta>0$ and $\mu>1$, then the origin of this system is globally
fixed-time stable and the settling time function $T$ can be estimated by $T\leq\frac{\pi\mu}{2\sqrt{\alpha\beta}}$.
\end{lemma}
\section{Design of Fixed-Time Frequency Estimator}
For the sinusoidal signal presented in (\ref{2.1}), we have
$y^{(1)}(t)=wBcos(wt+\phi_{0}),y^{(3)}(t)=-w^{3}Bcos(wt+\phi_{0})$, which will result in the following relation:
\begin{eqnarray}\label{3.1}
y^{(3)}(t)=-w^2y^{(1)}(t).
\end{eqnarray}
By integrating both side of (\ref{3.1}) over the time interval $[t-r,t]$ with design positive constant $r$ to be determined later, we have
\begin{eqnarray}\label{3.2}
\zeta\int_{t-r}^{t}|y^{(1)}(\tau)|\mathrm{d}\tau=\int_{t-r}^{t}|y^{(3)}(\tau)|\mathrm{d}\tau,\ \forall t\geq r,
\end{eqnarray}
where $\zeta=w^{2}$.
Defining two auxiliary variables
\begin{eqnarray}\label{3.3}
\gamma_{1}(t)=\int_{t-r}^{t}|y^{(1)}(\tau)|\mathrm{d}\tau,\ \gamma_{2}(t)=\int_{t-r}^{t}|y^{(3)}(\tau)|\mathrm{d}\tau,
\end{eqnarray}
then (\ref{3.2}) can be rewritten as
\begin{eqnarray}\label{3.4}
\zeta\gamma_{1}(t)=\gamma_{2}(t).
\end{eqnarray}
Note that the auxiliary variables $\gamma_{1}(t),\gamma_{2}(t)$ are unavailable for the unmeasurable signals $y^{(1)}(t)$ and $y^{(3)}(t)$. By substituting the estimation values of $y^{(1)}(t),y^{(3)}(t)$ into (\ref{3.3}), two new auxiliary variables $\hat{\gamma}_{1}(t),\hat{\gamma}_{2}(t)$ are obtained for $\forall t\geq0$ as:
\begin{eqnarray}\label{3.5}
\hat{\gamma}_{1}(t)=\int_{t-r}^{t}|z_2(\tau)|\mathrm{d}\tau,\ \hat{\gamma}_{2}(t)=\int_{t-r}^{t}|z_4(\tau)|\mathrm{d}\tau,
\end{eqnarray}
where $z_2(t)=z_4(t)=0$ for $\forall t\in[-r,0)$.

\subsection{Some Propositions}
In the following, some properties of the auxiliary variables $\hat{\gamma}_{1}(t),\hat{\gamma}_{2}(t)$ will be introduced firstly.

\begin{proposition}
The variables $\hat{\gamma}_{1}(t),\hat{\gamma}_{2}(t)$ are bounded and $\hat{\gamma}_{1}(t)=\gamma_{1}(t),\hat{\gamma}_{2}(t)=\gamma_{2}(t)$ holds for $\forall t\geq T_{1}+r$.
\end{proposition}

\emph{Proof of Proposition 1.}\ It can be concluded form Lemma 1 that $z_{i},i=1,\cdots,m$ are bounded and therefore $\hat{\gamma}_{1}(t),\hat{\gamma}_{2}(t)$ are also bounded. Note that
\begin{eqnarray}\label{3.6}
\hat{\gamma}_{1}(t)\!=\!\int_{t-r}^{t}|y^{(1)}(\tau)|\mathrm{d}\tau\!+\!\int_{t-r}^{t}\left(|z_2(\tau)|\!-\!|y^{(1)}(\tau)|\right)\mathrm{d}\tau.
\end{eqnarray}
Define $e_{\gamma_1}=\int_{t-r}^{t}\big||z_2(\tau)|-|y^{(1)}(\tau)|\big|\mathrm{d}\tau\geq0$. Then, we have
\begin{eqnarray}\label{3.7}
\gamma_{1}(t)-e_{\gamma_1}\leq\hat{\gamma}_{1}(t)\leq\gamma_{1}(t)+e_{\gamma_1}.
\end{eqnarray}
It is noted that when $t\geq T_{1}+r$, we have $e_{\gamma_1}=0$, which implies that
$\gamma_{1}(t)\leq\hat{\gamma}_{1}(t)\leq\gamma_{1}(t)$
holds for $\forall t\geq T_{1}+r$, i.e., $\hat{\gamma}_{1}(t)=\gamma_{1}(t)$ holds for $\forall t\geq T_{1}+r$. Similarly, we can also obtain that $\hat{\gamma}_{2}(t)=\gamma_{2}(t)$ holds for $\forall t\geq T_{1}+r$, which completes the proof of Proposition 1. \hfill $\blacksquare$


\begin{proposition}
When $w\neq0$ and $t\geq T_{1}+r$, the signal $\hat{\gamma}_{1}(t)$ satisfies a persistent excitation condition, i.e.,
for any given constant $r>0$, we can also find a constant $\epsilon>0$ such that $\hat{\gamma}_{1}(t)\geq\epsilon$
holds for $\forall t\geq T_{1}+r$.
\end{proposition}

\emph{Proof of Proposition 2.}\
When $t\geq T_{1}+r$, it follows from Proposition 1 that
\begin{eqnarray}\label{3.8.0}
\hat{\gamma}_{1}(t)=\gamma_{1}(t)
=Bw\int_{t-r}^{t}|cos(w\tau+\phi_{0})|\mathrm{d}\tau.
\end{eqnarray}
Note that the period of $|cos(wt+\phi_{0})|$ is $\pi/w$.
Thus, over any time interval $[t_{x1},t_{x2}]$ with $t_{x2}-t_{x1}\geq \pi k/w,k\in \mathbb{N}{+}$, the integral of $|cos(wt+\phi_{0})|$ satisfies
\begin{eqnarray}\label{3.8}
\!\int_{t_{x1}}^{t_{x2}}\!\!\!\!|cos(w\tau\!+\!\phi_{0})|\mathrm{d}\tau\!\geq\! k\!\!\int_{0}^{\pi/w}\!\!\!\!\!\!|cos(w\tau\!+\!\phi_{0})|\mathrm{d}\tau
\!=\!2k/w.
\end{eqnarray}
In the following, the conditions $r\geq\frac{\pi}{w}$ and $0<r<\frac{\pi}{w}$ will be discussed separately.

(i) When $r\geq\frac{\pi}{w}$, the time interval over $[t-r,t]$ satisfies $t-(t-r)=r\geq\frac{\pi}{w}$. Therefore, by using (\ref{3.8.0}) and (\ref{3.8}) with $k=1$, we have $\hat{\gamma}_{1}(t)\geq2B\geq2B_{min}$. Choosing $\epsilon=2B_{min}$, we have $\hat{\gamma}_{1}(t)\geq\epsilon$.

(ii) When $0<r<\frac{\pi}{w}$, one can calculate that
$\int_{t-r}^{t}|cos(w\tau+\phi_{0})|\mathrm{d}\tau\geq 2 \int_{0}^{\frac{r}{2}}|sinw\tau|\mathrm{d}\tau=\frac{2-2|cos\frac{rw}{2}|}{w}$. Along with (\ref{3.8.0}), we have $\hat{\gamma}_{1}(t)\geq2B(1-|cos\frac{rw}{2}|)\geq2B_{min}(1-|cos\frac{rw_{min}}{2}|)$. Letting $\epsilon=2B_{min}(1-|cos\frac{rw_{min}}{2}|)$, we have $\hat{\gamma}_{1}(t)\geq\epsilon$.

Hence, the proof of Proposition 2 is completed. \hfill $\blacksquare$
\subsection{Main Result}
Define the estimates $\hat{w},\hat{\zeta}(t)$ of $w,\zeta$, which is updated by the following adaptive law:
\begin{align}\label{3.9}
\dot{\hat{\zeta}}(t)&\!\!=\!\!\left\{
\begin{aligned}
&\!\!-\hat{\gamma}_{1}^{-1}(t)\bigg(\alpha_{1}e_{\gamma}^{1+\frac{q}{p}}+\beta_{1}e_{\gamma}^{1-\frac{q}{p}}+\hat{\zeta}(t)\big(|z_{2}(t)|\\
&\quad-\!\!|z_{2}(t\!-\!r)|\big)\!\!-\!\!\big(|z_{4}(t)|\!\!-\!\!|z_{4}(t\!-\!r)|\big)\bigg), \quad if\ \hat{\gamma}_{1}>\epsilon,\\
&\!\!-\alpha_{1}\hat{\zeta}^{1+\frac{q}{p}}(t)-\beta_{1}\hat{\zeta}^{1-\frac{q}{p}}(t),\quad\quad\ \quad\quad\quad otherwise, \\
\end{aligned}
\right.\nonumber\\
\hat{w}(t)&\!\!=\!\!|\hat{\zeta}(t)|^{1/2},
\end{align}
where
\begin{eqnarray}\label{3.9.0}
e_{\gamma}=\hat{\zeta}(t)\hat{\gamma}_{1}(t)-\hat{\gamma}_{2}(t),
\end{eqnarray}
and $\alpha_{1},\beta_{1}$ are positive constant, $0<q<2p$ are odd integers, $\epsilon>0$ is selected according to Proposition 2, $\hat{\gamma}_{1}(t)$ is defined in (\ref{3.5}), $z_{2}(t),z_{4}(t)$ are the states of the observer (\ref{2.2}).
The frequency estimator can be implemented as (\ref{2.2}), (\ref{3.5}), (\ref{3.9}) which will result in the following theorem.

\begin{theorem}
For the biased sinusoidal signal $y(t)$ defined in (\ref{2.1}), if Assumption 1 is satisfied, then the frequency estimator (\ref{2.2}), (\ref{3.5}), (\ref{3.9}) can achieve fixed-time frequency estimation, i.e., there exists a time $T_{max}$ independent of initial condition such that $\hat{w}(t)=w$ holds for $\forall t\geq T_{max}$.
\end{theorem}


\emph{Proof of Theorem 1.}\
For the condition $w\in[w_{min},+\varpropto)$, the derivative of $e_{\gamma}$ can be calculated with (\ref{3.5}) as
\begin{eqnarray}\label{3.10}
\dot{e}_{\gamma}&=&\dot{\hat{\zeta}}(t)\hat{\gamma}_{1}(t)+\hat{\zeta}(t)\dot{\hat{\gamma}}_{1}(t)-\dot{\hat{\gamma}}_{2}(t)\nonumber\\
&=&\dot{\hat{\zeta}}(t)\hat{\gamma}_{1}(t)+\hat{\zeta}(t)(|z_{2}(t)|-|z_{2}(t-r)|)\nonumber\\
&&-(|z_{4}(t)|-|z_{4}(t-r)|).
\end{eqnarray}
It follows from Proposition 2 that $\hat{\gamma}_{1}(t)\geq\epsilon$ holds for $\forall t\geq T_{1}+r$. Therefore, when $t\geq T_{1}+r$, substituting the adaptive law (\ref{3.9}) into (\ref{3.10}) will lead to
\begin{eqnarray}\label{3.11}
\dot{e}_{\gamma}=-\alpha_{1}e_{\gamma}^{1+\frac{q}{p}}-\beta_{1}e_{\gamma}^{1-\frac{q}{p}}.
\end{eqnarray}
The derivative of Lyapunov function $V=\frac{1}{2}e_{\gamma}^2$ along (\ref{3.11}) satisfies
\begin{eqnarray}\label{3.12}
\dot{V}=-2^{1+\frac{q}{2p}}\alpha_{1}V^{1+\frac{q}{2p}}-2^{1-\frac{q}{2p}}\beta_{1}V^{1-\frac{q}{2p}},
\end{eqnarray}
which implies that $V$ and thus $e_{\gamma}$ will converge to zero in a time $T_{max2}$ independent of initial condition, i.e., $e_{\gamma}=0$ holds for $\forall t\geq T_{1}+T_{max2}+r$.
According to (\ref{3.9.0}), $e_{\gamma}=0$ means $\hat{\zeta}(t)\hat{\gamma}_{1}(t)=\hat{\gamma}_{2}(t)$. Note that  $\hat{\gamma}_{1}(t)=\gamma_{1}(t),\gamma_{2}(t)=\hat{\gamma}_{2}(t)$ holds for $\forall t\geq T_{1}+r$. Therefore, when $t\geq T_{1}+T_{max2}+r$, we have
\begin{eqnarray}\label{3.13}
\hat{\zeta}(t)\gamma_{1}(t)=\gamma_{2}(t).
\end{eqnarray}
Subtracting (\ref{3.13}) with (\ref{3.4}), we have $(\zeta-\hat{\zeta}(t))\gamma_{1}(t)=0$, which means that $\hat{\zeta}(t)=\zeta$ and thus $\hat{w}(t)=w$ holds for $\forall t\geq T_{1}+T_{max2}+r$.

When $w=0$, it is easy to verify that $\hat{\gamma}_{1}(t)=\gamma_{1}(t)=0$ holds for $\forall t\geq T_{1}+r$. Then, when $t\geq T_{1}+r$, the adaptive law (\ref{3.9}) reduces to
$\dot{\hat{\zeta}}(t)=-\alpha_{1}\hat{\zeta}^{1+\frac{q}{p}}(t)-\beta_{1}\hat{\zeta}^{1-\frac{q}{p}}(t)$, which is fixed-time stable ({\color{blue}{Polyakov, 2012}}). Similar to previous analysis, we can conclude that there exists a fixed-time $T_{max1}$ such that $\hat{\zeta}(t)=0$ holds for $\forall t\geq T_{1}+T_{max1}+r$. Note that $\hat{\zeta}(t)=0$ implies $\hat{w}(t)=w=0$ holds for $\forall t\geq T_{1}+T_{max1}+r$.

Define $T_{max}=T_{1}+r+\max\{T_{max1},T_{max2}\}$. We can conclude that $\hat{w}=w$ will be established for $\forall t\geq T_{max}$, which completes the proof of Theorem 1. \hfill $\blacksquare$

\section{Robustness Analysis}
In practice, one never has access to perfect measurements. Therefore, the robustness analysis of the proposed algorithm in the presence of measurement noise will be given in the following.
Suppose that the signal $y(t)$ is measured in the presence of bounded measurement noise $n(t):\ |n(t)|\leq \eta$, i.e., the measurement $\hat{y}(t)$ satisfies $\hat{y}(t)=y(t)+n(t)$. By replacing $y(t)$ with the $\hat{y}(t)$ in the observer (\ref{2.2}), the following accuracies
\begin{align}\label{4.1}
\!\!|\tilde{z}_{i}|\!\!\triangleq\!\!|z_{i}\!\!-\!\!y^{(i-1)}(t)|\!\leq \! \mathcal{O}(\eta^{(m-i+1)/m}),i=1,\cdots,m,
\end{align}
can be obtained after finite-time ({\color{blue}{Angulo, et al., 2013}}).
In view of (\ref{4.1}), the following accuracies
\begin{align}\label{4.2}
|\tilde{\gamma}_{1}(t)|\triangleq|\gamma_{1}(t)-\hat{\gamma}_{1}(t)|\leq \mathcal{O}(\eta^{(m-1)/m}),\nonumber\\
|\tilde{\gamma}_{2}(t)|\triangleq|\gamma_{2}(t)-\hat{\gamma}_{2}(t)|\leq \mathcal{O}(\eta^{(m-3)/m}),
\end{align}
can also be obtained after finite-time.

Note that for the existence of noise, the persistent excitation condition in Proposition 2 may not be satisfied. This make the estimator (\ref{3.9}) inactive or only active for some instants, which is separated by the threshold $\epsilon$. Take a very particular case for example, $n(t)=-y(t)$ and $\hat{y}(t)=0$. In this case, the perturbed measurement $\hat{y}(t)$ cannot be used to estimate the frequency. Latter, we will show that when the the persistent excitation condition is still satisfied even in the presence of noise, the proposed frequency estimator is ISS with respect the measurement noise $n(t)$.

It follows from the proof of Theorem 1 that $e_{\gamma}=0$ and thus $\hat{\zeta}(t)\hat{\gamma}_{1}(t)=\hat{\gamma}_{2}(t)$ can be established after finite time. Equation (\ref{3.4}) can be rewritten as $\zeta\hat{\gamma}_{1}(t)-\zeta\tilde{\gamma}_{1}(t)=\hat{\gamma}_{2}(t)-\tilde{\gamma}_{2}(t)$. Then, subtracting the above two equations, we have
$|\tilde{\zeta}(t)|\triangleq|\hat{\zeta}(t)-\zeta|
\leq|\hat{\gamma}_{1}^{-1}(t)(\tilde{\gamma}_{2}(t)-\zeta\tilde{\gamma}_{1}(t))|
\leq\epsilon^{-1}|\tilde{\gamma}_{2}(t)|+|\zeta||\tilde{\gamma}_{1}(t)|
$. Therefore, $\tilde{\zeta}(t)$ and the estimator error $\tilde{w}(t)=\hat{w}(t)-w$ is ISS with respect the measurement noise $n(t)$, which can be summarized as follows:
\begin{corollary}
Suppose that the measurement $\hat{y}(t)$ still satisfies the persistent excitation condition in Proposition 2. Then, the estimator error $\tilde{w}(t)=\hat{w}(t)-w$ is ISS with respect to any bounded measurement noise $n(t)$.
\end{corollary}

\section{Simulation}
\begin{figure}[htbp]
\centering
\subfloat{
\begin{minipage}[t]{0.5\textwidth}
\centering
\rotatebox{360}{\scalebox{0.5}[0.5]{\includegraphics{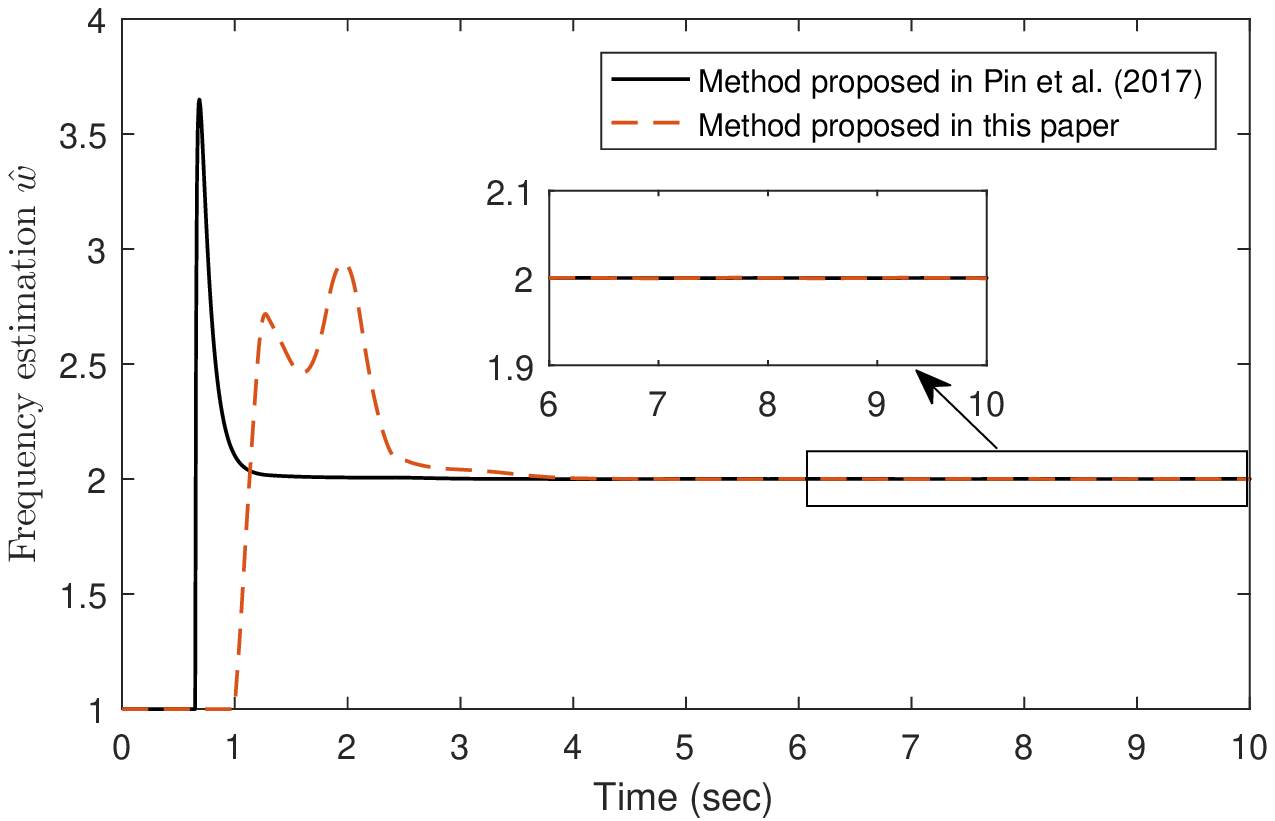}}}
\end{minipage}
}\\
\subfloat{
\begin{minipage}[t]{0.5\textwidth}
\centering
\rotatebox{360}{\scalebox{0.5}[0.5]{\includegraphics{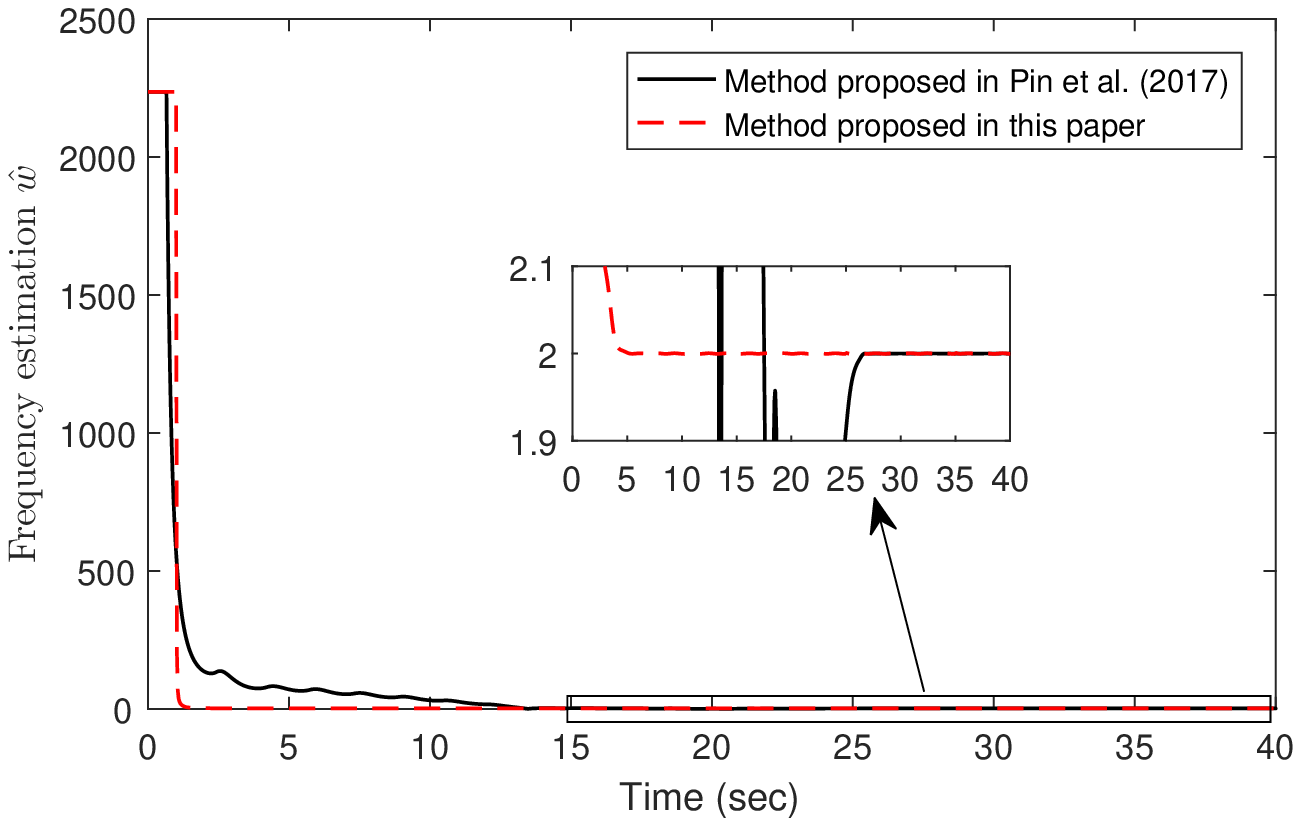}}}
\end{minipage}
}
\center{\small Fig. 1.\ Frequency estimation in the absence of noise by using the proposed method and the method proposed in {\color{blue}{Pin et al. (2017)}}: Above. small initial estimation error condition $\tilde{w}(0)=-1$; Below. large initial estimation error condition $\tilde{w}(0)\approx2235$.}
\end{figure}

\begin{figure}[htbp]
\centering
\subfloat{
\begin{minipage}[t]{0.5\textwidth}
\centering
\rotatebox{360}{\scalebox{0.5}[0.5]{\includegraphics{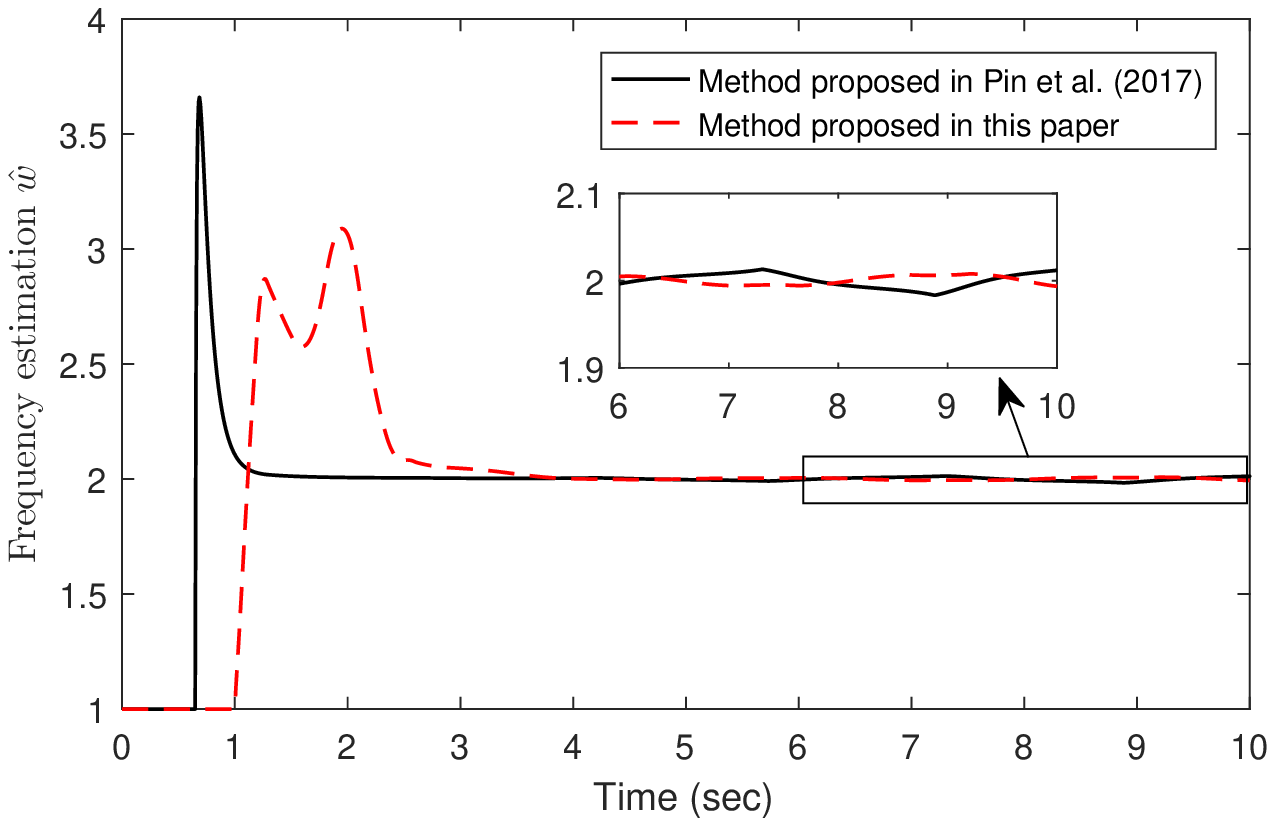}}}
\end{minipage}
}\\
\subfloat{
\begin{minipage}[t]{0.5\textwidth}
\centering
\rotatebox{360}{\scalebox{0.5}[0.5]{\includegraphics{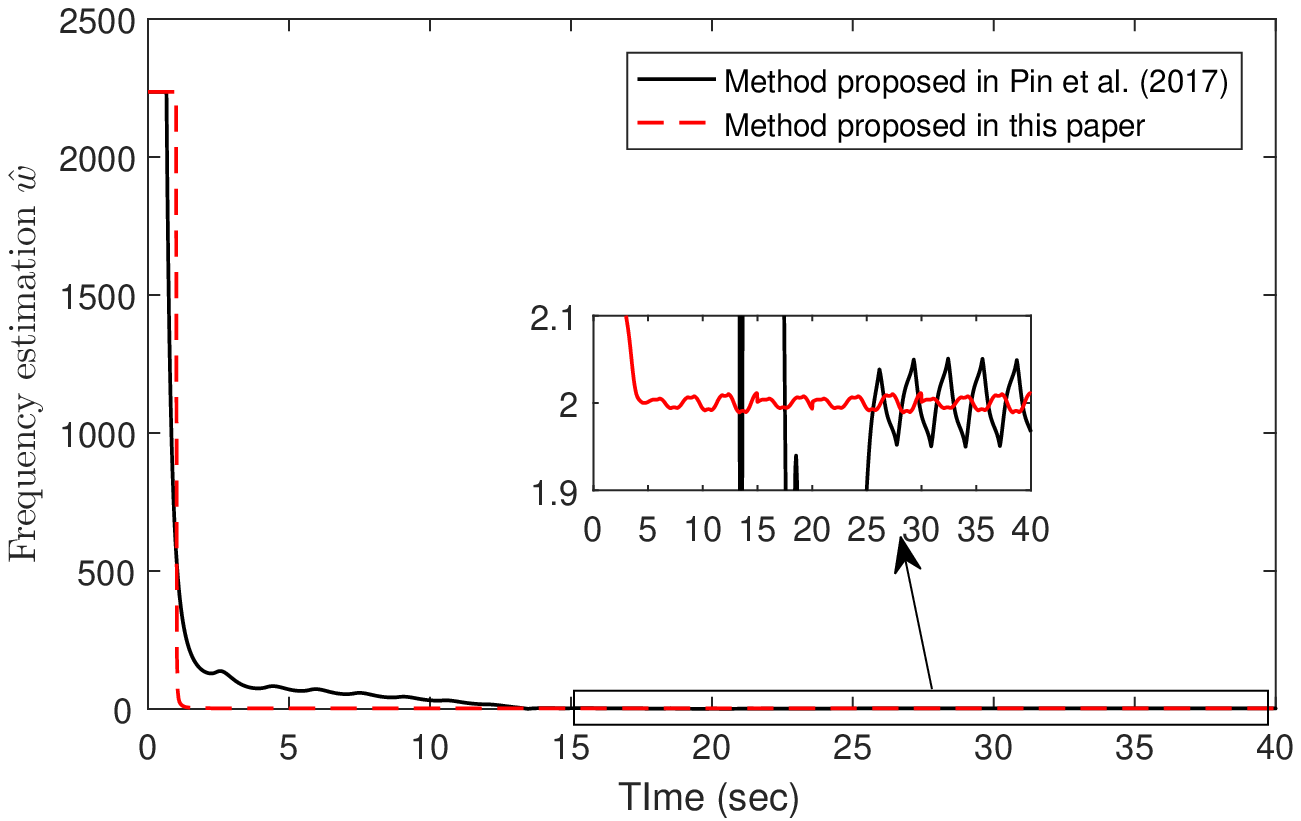}}}
\end{minipage}
}
\center{\small Fig. 2.\ Frequency estimation in the presence of noise by using the proposed method and the method proposed in {\color{blue}{Pin et al. (2017)}}: Above. small initial estimation error condition $\tilde{w}(0)=-1$; Below. large initial estimation error condition $\tilde{w}(0)\approx2235$.}
\end{figure}
In this section, simulation results for the frequency estimation of the signal $y(t)=4sin(2t+2)+10$ will be given.
For the observer (\ref{2.2}), we select $\kappa_{1}=16,\kappa_{2}=88,\kappa_{3}=140,\kappa_{4}=110$, $k_{1}=24,k_{2}=216,k_{3}=864,k_{4}=1296$, $T_{u}=3$ and $\alpha=0.6$. Select $r=1$ for the auxiliary variables $\hat{\gamma}_{1}(t),\hat{\gamma}_{2}(t)$ defined in (\ref{3.5}). For the adaptive law (\ref{3.9}), we select $p=3,q=1,\epsilon=0.01$ and $\alpha_{1}=\beta_{1}=1$.
To show the effectiveness of the proposed method, the finite-time adaptive frequency estimator proposed in {\color{blue}{Pin et al. (2017)}} will be used to make simulation comparison with parameters selected as $\beta_1=1,\beta_2=2,\beta_3=3,\bar{\beta}=2.5,g=0.1,g_a=25,L_1=1.5,L_2=1.1$ and $\delta_{\epsilon}=0.001$.

Firstly, simulation results in the absence of measurement noise by using the proposed method and the method proposed in {\color{blue}{Pin et al. (2017)}} is given in Fig. 1. Different initial conditions are used to make comparison. It can be clearly observed from Fig. 1 that
the proposed method can achieve exact estimation of the frequency within a fixed-time $5s$ no matter how large the initial values are selected, while the settling time of the method in {\color{blue}{Pin et al. (2017)}} grows from $5s$ to $25s$ when the initial condition grow.
To show the robustness of the proposed method, a bounded measurement noise $|n(t)|\leq0.25$ is considered in Fig. 2. It can be observed from Fig. 2 that similar to the existing method, our proposed method is also ISS with respect to bounded measurement noise $n(t)$.

\textbf{Note:} More details about the simulation can be found in Appendix attached at the end of the manuscript.
\section{Conclusion}
This technical communique has developed a fixed-time frequency estimator for biased sinusoidal signals for the first time. How to extend the result to handle multiple biased and damped sinusoidal signals is the future work.


\clearpage
\section*{Appendix A. M-files and Simulation Results}
The detailed m-files of the simulation results presented in Figs. 1-2 are given in the Appendix A.1 and A.2, respectively. Moreover, the m-files and simulation results for Example 3 in {\color{blue}{Pin et al. (2017)}} with different initial values $\tilde{w}(0)$ are also given in Appendix A.3. It can observed from Fig. A.1 in Appendix A.3  that the settling time of the method in {\color{blue}{Pin et al. (2017)}} grows from $1s$ to more than $700s$ when the initial error conditions grow.

\emph{\textbf{Appendix A.1: M-file of the proposed algorithm to generate the simulation result in this paper}}
\begin{lstlisting}[style=Matlab-editor]
clear all;
close all;
clc

dt=0.000001;T=10;L=160;A=10;B=4;
w=2;r=1;afa1=1;beta1=1;Tu=1;afa=0.6;
q1=1;p1=3;q2=1;p2=3;g1=0;g2=0;
ka1=5*L^(1/4);ka2=3*L^(1/3)*ka1^(2/3);
ka3=1.5*L^(1/2)*ka2^(1/2);ka4=1.1*L;
k1=24;k2=216;k3=864;k4=1296;
z1(1)=0;z2(1)=0;z3(1)=0;z4(1)=0;

hy=1;%Small initial error condition
%hy=5*10^6;%Large initial error condition

for i=1:T/dt
   t=i*dt;
   theta=(sign(t-Tu)+1)/2;

   y(1)=0; y(i+1)=B*sin(w*t+2)+A;
   %hhy(i)=B*sin(w*t+2)+A+2*cos(0.05*t);
   y1(1)=0; y1(i+1)=B*w*cos(w*t+2);
   y2(1)=0; y2(i+1)=-B*w^2*sin(w*t+2);
   y3(1)=0; y3(i+1)=-B*w^3*cos(w*t+2);

   wz1(i)=z1(i)-y(i);
   %wz1(i)=z1(i)-hhy(i);
   z1(i+1)=(-ka1*theta*sign(wz1(i))*(abs(wz1(i)))^(3/4)-k1*(1-theta)*sign(wz1(i))*(abs(wz1(i)))^((4+afa)/4)+z2(i))*dt+z1(i);
   z2(i+1)=(-ka2*theta*sign(wz1(i))*(abs(wz1(i)))^(2/4)-k2*(1-theta)*sign(wz1(i))*(abs(wz1(i)))^((4+2*afa)/4)+z3(i))*dt+z2(i);
   z3(i+1)=(-ka3*theta*sign(wz1(i))*(abs(wz1(i)))^(1/4)-k3*(1-theta)*sign(wz1(i))*(abs(wz1(i)))^((4+3*afa)/4)+z4(i))*dt+z3(i);
   z4(i+1)=(-ka4*theta*sign(wz1(i))-k4*(1-theta)*sign(wz1(i))*(abs(wz1(i)))^((4+4*afa)/4))*dt+z4(i);

   if i<r/dt+1

    g1=abs(z2(i))*dt+g1;
    g2=abs(z4(i))*dt+g2;
    hy=hy;

   else

    h=i-1000000;
    a(i)=z2(i-1000000);
    g1=(abs(z2(i))-abs(z2(i-1000000)))*dt+g1;
    g2=(abs(z4(i))-abs(z4(i-1000000)))*dt+g2;

    e(i)=hy*g1-g2;

    if g1>0.01

     v=-1/g1*(afa1*e(i)^(1+q2/p2)+beta1*e(i)^(1-q1/p2)+hy*(abs(z2(i))-abs(z2(i-1000000)))-(abs(z4(i))-abs(z4(i-1000000))));
     hy=v*dt+hy;

    else

     hy=(-afa1*hy^(1+q2/p2)-beta1*hy^(1-q1/p1))*dt+hy;

    end

   end

   hw(i)=abs(hy)^(1/2);
   w1(i)=w;

end
t=dt:dt:T;
figure;
plot(t,hw,t,w1);
\end{lstlisting}

\emph{\textbf{Appendix A.2: M-file of {\color{blue}{Pin et al. (2017)}} to generate the simulation result in this paper }}
\begin{lstlisting}[style=Matlab-editor]
clear all;
close all;
clc

dt=0.000001;T=40;A=10;B=4;w=2;
b1=1;b2=2;b3=3;b=2.5;g=1;deta=0.001;L1=10;L2=2;
xi11(1)=0;xi13(1)=0;xi21(1)=0;xi23(1)=0;xi31(1)=0;xi33(1)=0;
r1(1)=0;r2(1)=0;yo(1)=0;

%ho(1)=1;%small initial error condition
ho(1)=5*10^6;%large initial error condition

for i=1:T/dt
   t=i*dt;
   %signal y=B*sin(w*t+2)+A
   % y(i)=B*sin(w*t+2)+A+2*cos(0.05*t);
   y(i)=B*sin(w*t+2)+A;

   %signal F1(t,t) and its derivatives F11,F12,F13 (in view of Eq.(16))
   F1=1-3*exp(-b*t)+3*exp(-2*b*t)-exp(-3*b*t);
   F11=b1-3*(b1-b)*exp(-b*t)+3*(b1-2*b)*exp(-2*b*t)-(b1-3*b)*exp(-3*b*t);
   F12=b1^2-3*(b1-b)^2*exp(-b*t)+3*(b1-2*b)^2*exp(-2*b*t)-(b1-3*b)^2*exp(-3*b*t);
   F13=b1^3-3*(b1-b)^3*exp(-b*t)+3*(b1-2*b)^3*exp(-2*b*t)-(b1-3*b)^3*exp(-3*b*t);

   %signal F2(t,t) and its derivatives F21,F22,F23 (in view of Eqs.(16) and (20))
   F2=F1;
   F21=b2-3*(b2-b)*exp(-b*t)+3*(b2-2*b)*exp(-2*b*t)-(b2-3*b)*exp(-3*b*t);
   F22=b2^2-3*(b2-b)^2*exp(-b*t)+3*(b2-2*b)^2*exp(-2*b*t)-(b2-3*b)^2*exp(-3*b*t);
   F23=b2^3-3*(b2-b)^3*exp(-b*t)+3*(b2-2*b)^3*exp(-2*b*t)-(b2-3*b)^3*exp(-3*b*t);

   %signal F3(t,t) and its derivatives F31,F32,F33 (in view of Eqs.(16) and (20))
   F3=F1;
   F31=b3-3*(b3-b)*exp(-b*t)+3*(b3-2*b)*exp(-2*b*t)-(b3-3*b)*exp(-3*b*t);
   F32=b3^2-3*(b3-b)^2*exp(-b*t)+3*(b3-2*b)^2*exp(-2*b*t)-(b3-3*b)^2*exp(-3*b*t);
   F33=b3^3-3*(b3-b)^3*exp(-b*t)+3*(b3-2*b)^3*exp(-2*b*t)-(b3-3*b)^3*exp(-3*b*t);

   %auxiliary systems (22) to generate auxiliary signals in (19)
   xi11(i+1)=(F11*y(i)-b1*xi11(i))*dt+xi11(i);
   xi13(i+1)=(F13*y(i)-b1*xi13(i))*dt+xi13(i);
   xi21(i+1)=(F21*y(i)-b2*xi21(i))*dt+xi21(i);
   xi23(i+1)=(F23*y(i)-b2*xi23(i))*dt+xi23(i);
   xi31(i+1)=(F31*y(i)-b3*xi31(i))*dt+xi31(i);
   xi33(i+1)=(F33*y(i)-b3*xi33(i))*dt+xi33(i);

   %auxiliary signals K (in view of Eq. (19))
   K1a=xi13(i)-F12*y(i);K2a=xi23(i)-F22*y(i);K3a=xi33(i)-F32*y(i);
   K1b=F11;K2b=F21;K3b=F31;
   K1d=xi11(i)-F1*y(i);K2d=xi21(i)-F2*y(i);K3d=xi31(i)-F3*y(i);

   %vector form (inview of Eq. (26))
   Ka=[K1a;K2a;K3a];Kd=[K1d;K2d;K3d];F=[K3b-K2b;K1b-K3b;K2b-K1b];
   K1=Ka'*F;K2=Kd'*F;

   %Deformation of the system (In view of (27)-(29))
   dr1=abs(K1)-g*r1(i);dr2=abs(K2)-g*r2(i);
   r1(i+1)=dr1*dt+r1(i);r2(i+1)=dr2*dt+r2(i);

   %adaptive law (In view of (32))
   Ro(i)=r1(i)-r2(i)*ho(i);
   if r2(i)>deta
       ho(i+1)=(yo(i)+L1*abs(Ro(i))^(1/2)*sign(Ro(i))-ho(i)*dr2+dr1)/r2(i)*dt+ho(i);
   else
       ho(i+1)=ho(i);
   end
   yo(i+1)=(L2*sign(Ro(i)))*dt+yo(i);

%Frequency estimation
hw(i)=abs(ho(i))^(1/2);
w1(i)=w;

end
t=dt:dt:T;
figure;
plot(t,hw,t,w1);
\end{lstlisting}

\emph{\textbf{Appendix A.3: M-file and Simulation Results for Example 3 in {\color{blue}{Pin et al. (2017)}}}}
\begin{lstlisting}[style=Matlab-editor]
clear all;
close all;
clc

dt=0.000001;T=2;A=2;B=3;w=4;b1=1;b2=2;
b3=3;b=2.5;g=3;deta=0.0001;L1=30;L2=2;
xi11(1)=0;xi13(1)=0;xi21(1)=0;
xi23(1)=0;xi31(1)=0;xi33(1)=0;
r1(1)=0;r2(1)=0;yo(1)=0;

ho(1)=5^(1/2);%Small initial error condition
%ho(1)=5*10^2;%Initial error condition hw-w=17.3607
%ho(1)=5*10^9;%Initial error condition hw-w=7*10^4
%ho(1)=5*10^10;%Initial error condition hw-w=2.2*10^5
%ho(1)=5*10^11;%Initial error condition hw-w=7.07*10^5

for i=1:T/dt

   t=i*dt;
   y(i)=B*sin(w*t+pi/4)+A;

   %signal F1(t,t) and its derivatives F11,F12,F13 (in view of Eq.(16))
   F1=1-3*exp(-b*t)+3*exp(-2*b*t)-exp(-3*b*t);
   F11=b1-3*(b1-b)*exp(-b*t)+3*(b1-2*b)*exp(-2*b*t)-(b1-3*b)*exp(-3*b*t);
   F12=b1^2-3*(b1-b)^2*exp(-b*t)+3*(b1-2*b)^2*exp(-2*b*t)-(b1-3*b)^2*exp(-3*b*t);
   F13=b1^3-3*(b1-b)^3*exp(-b*t)+3*(b1-2*b)^3*exp(-2*b*t)-(b1-3*b)^3*exp(-3*b*t);

   %signal F2(t,t) and its derivatives F21,F22,F23 (in view of Eqs.(16) and (20))
   F2=F1;
   F21=b2-3*(b2-b)*exp(-b*t)+3*(b2-2*b)*exp(-2*b*t)-(b2-3*b)*exp(-3*b*t);
   F22=b2^2-3*(b2-b)^2*exp(-b*t)+3*(b2-2*b)^2*exp(-2*b*t)-(b2-3*b)^2*exp(-3*b*t);
   F23=b2^3-3*(b2-b)^3*exp(-b*t)+3*(b2-2*b)^3*exp(-2*b*t)-(b2-3*b)^3*exp(-3*b*t);

   %signal F3(t,t) and its derivatives F31,F32,F33 (in view of Eqs.(16) and (20))
   F3=F1;
   F31=b3-3*(b3-b)*exp(-b*t)+3*(b3-2*b)*exp(-2*b*t)-(b3-3*b)*exp(-3*b*t);
   F32=b3^2-3*(b3-b)^2*exp(-b*t)+3*(b3-2*b)^2*exp(-2*b*t)-(b3-3*b)^2*exp(-3*b*t);
   F33=b3^3-3*(b3-b)^3*exp(-b*t)+3*(b3-2*b)^3*exp(-2*b*t)-(b3-3*b)^3*exp(-3*b*t);

   %auxiliary systems (22) to generate auxiliary signals in (19)
   xi11(i+1)=(F11*y(i)-b1*xi11(i))*dt+xi11(i);
   xi13(i+1)=(F13*y(i)-b1*xi13(i))*dt+xi13(i);
   xi21(i+1)=(F21*y(i)-b2*xi21(i))*dt+xi21(i);
   xi23(i+1)=(F23*y(i)-b2*xi23(i))*dt+xi23(i);
   xi31(i+1)=(F31*y(i)-b3*xi31(i))*dt+xi31(i);
   xi33(i+1)=(F33*y(i)-b3*xi33(i))*dt+xi33(i);

   %auxiliary signals K (in view of Eq. (19))
   K1a=xi13(i)-F12*y(i);K2a=xi23(i)-F22*y(i);K3a=xi33(i)-F32*y(i);
   K1b=F11;K2b=F21;K3b=F31;
   K1d=xi11(i)-F1*y(i);K2d=xi21(i)-F2*y(i);K3d=xi31(i)-F3*y(i);

   %vector form (inview of Eq. (26))
   Ka=[K1a;K2a;K3a];Kd=[K1d;K2d;K3d];F=[K3b-K2b;K1b-K3b;K2b-K1b];
   K1=Ka'*F;K2=Kd'*F;

   %Deformation of the system (In view of (27)-(29))
   dr1=abs(K1)-g*r1(i);dr2=abs(K2)-g*r2(i);
   r1(i+1)=dr1*dt+r1(i);r2(i+1)=dr2*dt+r2(i);

   %adaptive law (In view of (32))
   Ro(i)=r1(i)-r2(i)*ho(i);

   if r2(i)>deta
       ho(i+1)=(yo(i)+L1*abs(Ro(i))^(1/2)*sign(Ro(i))-ho(i)*dr2+dr1)/r2(i)*dt+ho(i);
   else
       ho(i+1)=ho(i);
   end
   yo(i+1)=(L2*sign(Ro(i)))*dt+yo(i);

 %Frequency estimation
hw(i)=abs(ho(i))^(1/2);
w1(i)=w;

end
t=dt:dt:T;
figure;
plot(t,hw,t,w1);
\end{lstlisting}

\begin{figure*}[htbp]
\centering
\subfloat{
\begin{minipage}[t]{0.5\textwidth}
\centering
\rotatebox{360}{\scalebox{0.5}[0.5]{\includegraphics{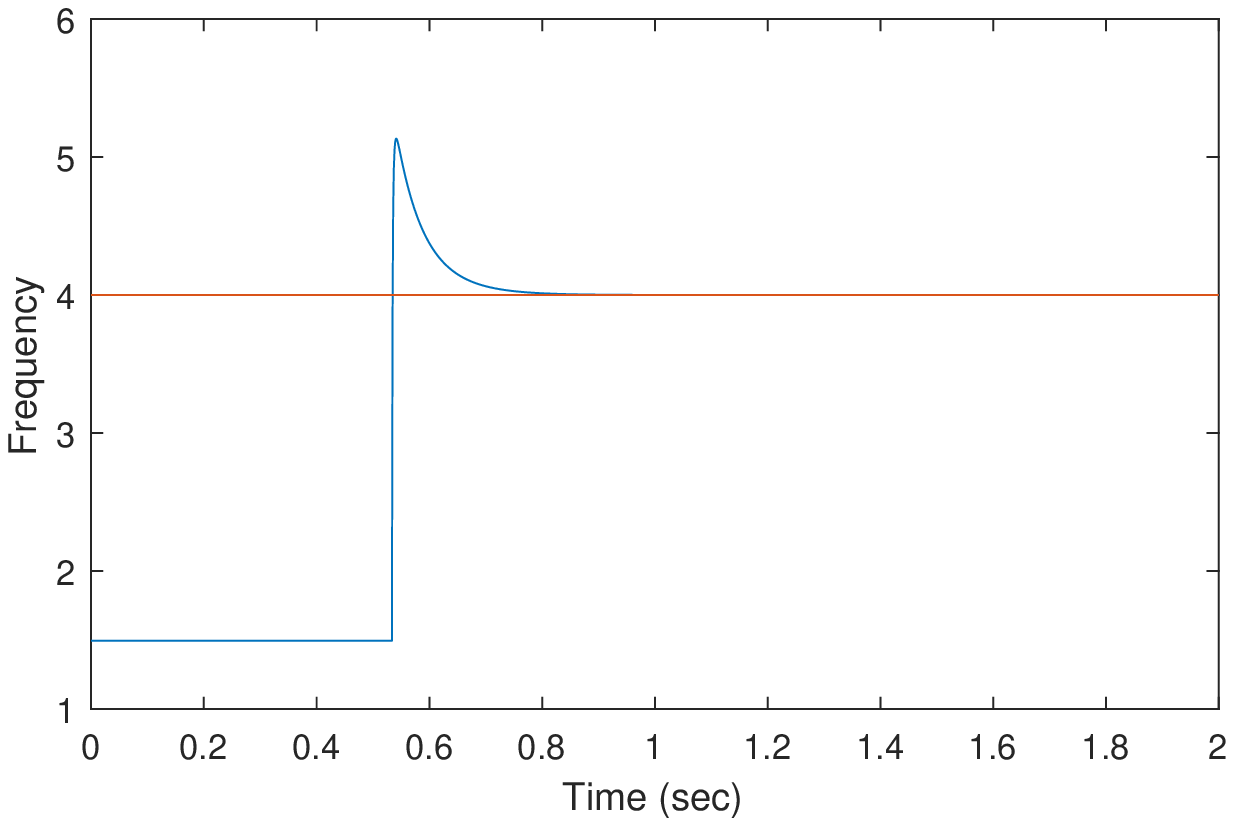}}}
\center{\small Scenario $(a)$: $\tilde{w}(0)\approx-1.7$.}
\end{minipage}
}
\subfloat{
\begin{minipage}[t]{0.5\textwidth}
\centering
\rotatebox{360}{\scalebox{0.5}[0.5]{\includegraphics{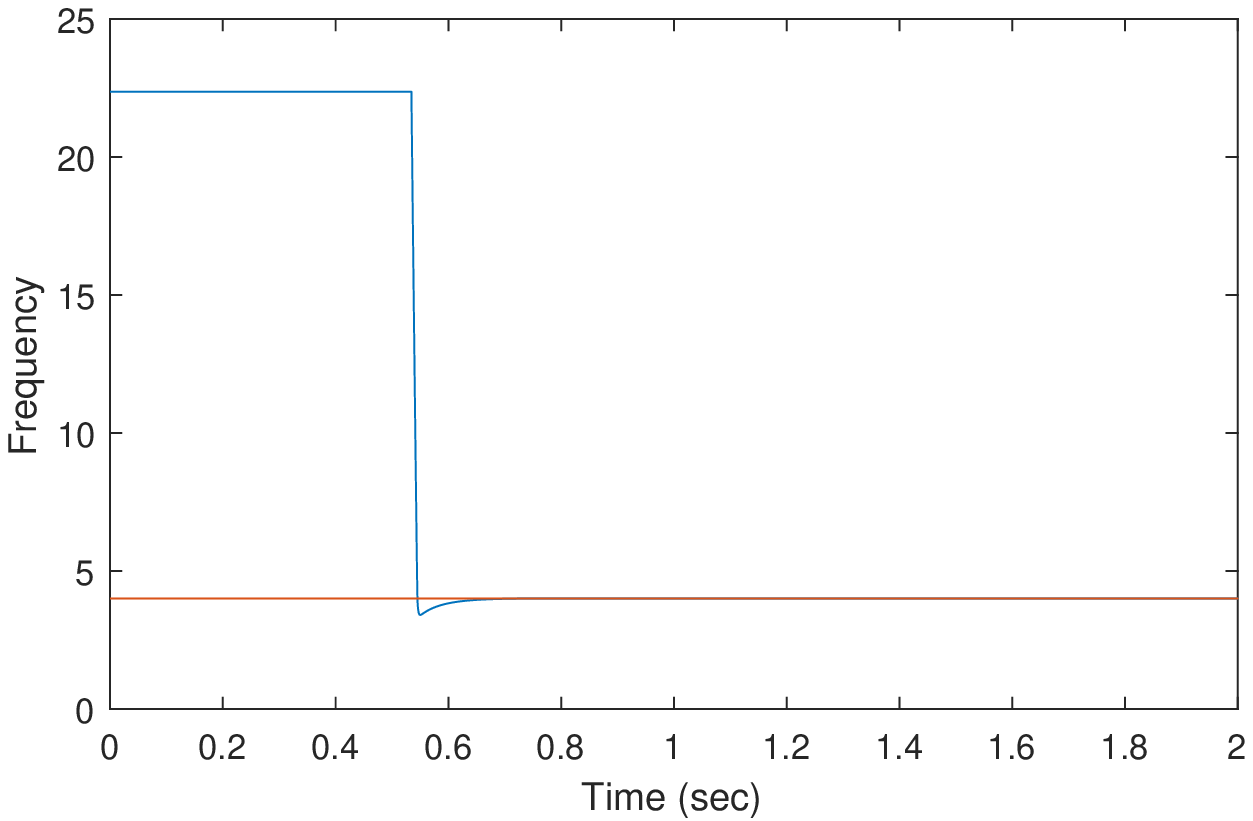}}}
\center{\small Scenario $(b)$: $\tilde{w}(0)\approx17$.}
\end{minipage}
}\\
\subfloat{
\begin{minipage}[t]{0.5\textwidth}
\centering
\rotatebox{360}{\scalebox{0.5}[0.5]{\includegraphics{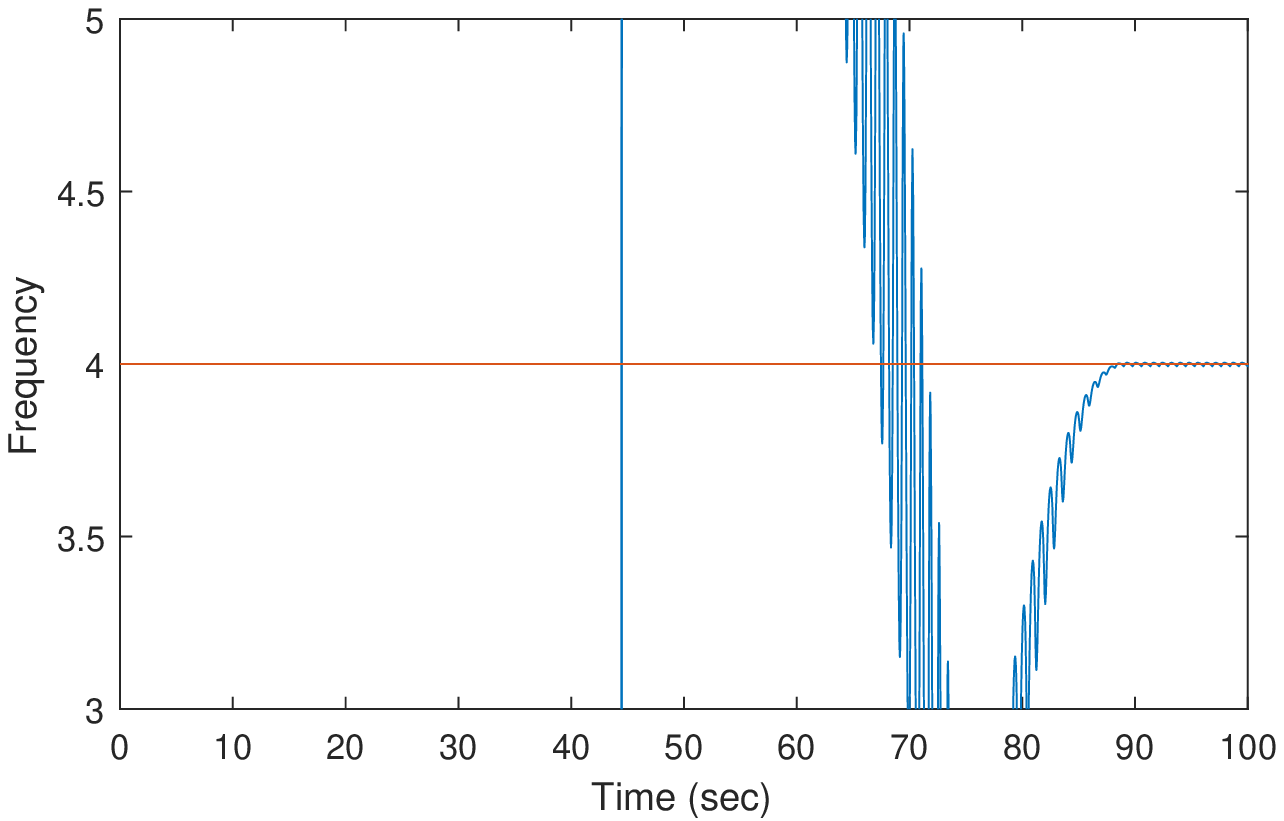}}}
\center{\small Scenario $(c)$: $\tilde{w}(0)\approx7*10^4$.}
\end{minipage}
}
\subfloat{
\begin{minipage}[t]{0.5\textwidth}
\centering
\rotatebox{360}{\scalebox{0.5}[0.5]{\includegraphics{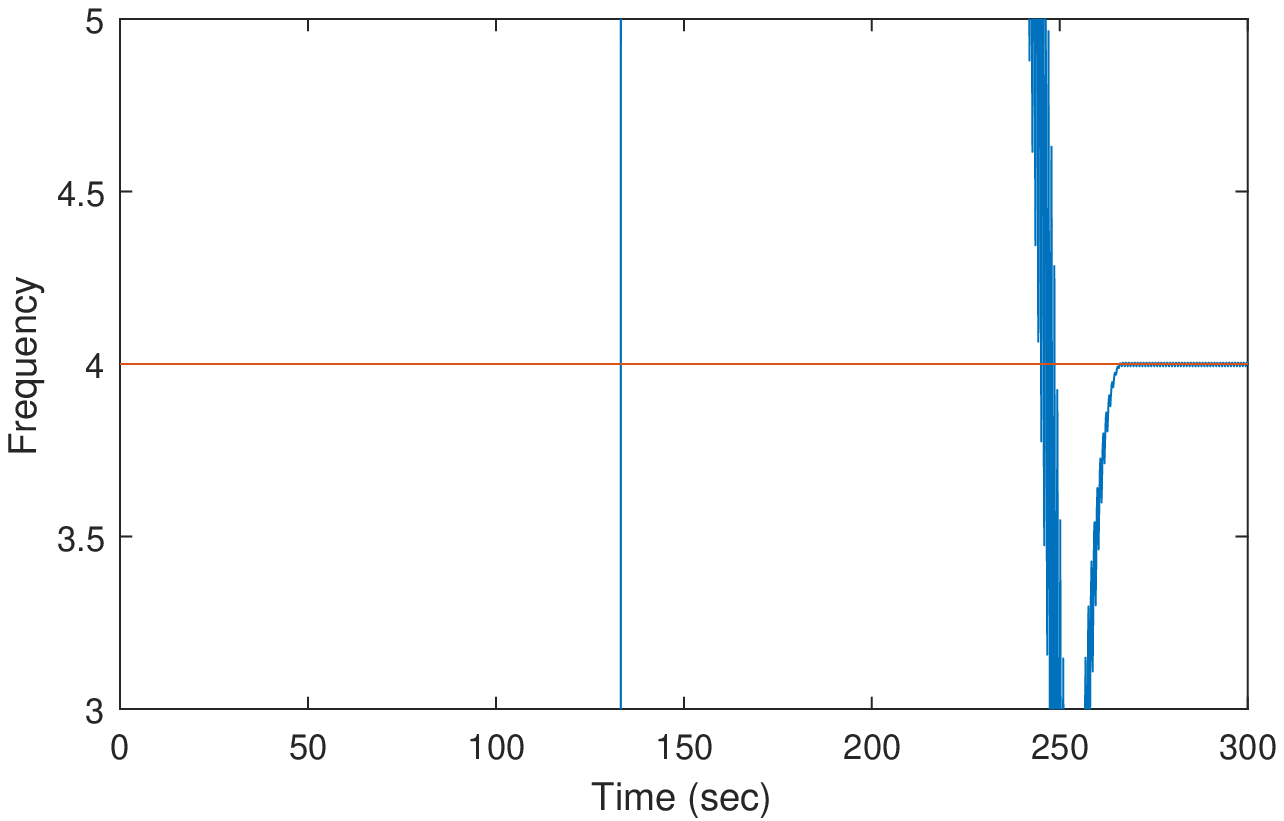}}}
\center{\small Scenario $(d)$: $\tilde{w}(0)\approx2.2*10^5$.}
\end{minipage}
}\\
\subfloat{
\begin{minipage}[t]{0.5\textwidth}
\centering
\rotatebox{360}{\scalebox{0.5}[0.5]{\includegraphics{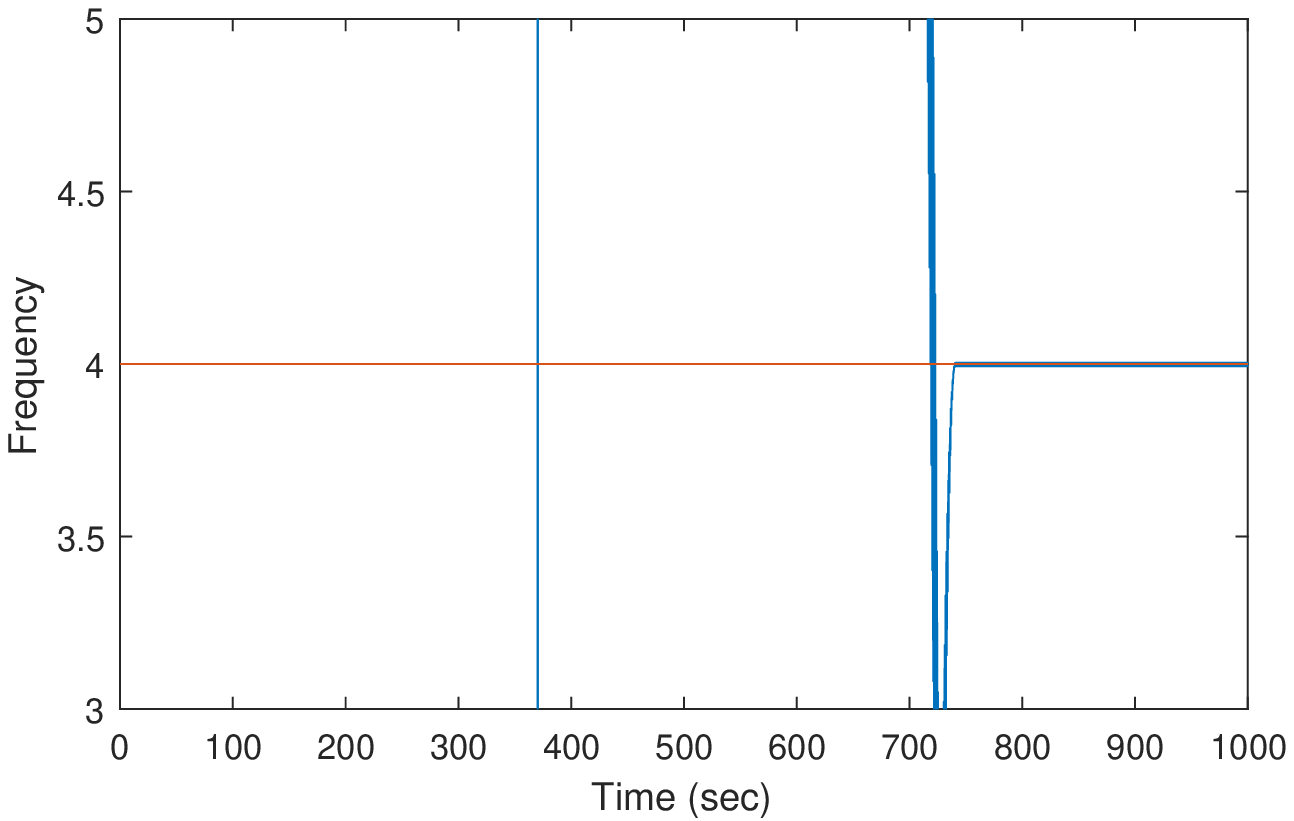}}}
\center{\small Scenario $(e)$: $\tilde{w}(0)\approx7*10^6$.}
\end{minipage}
}
\center{\small Fig. A.1.\ Simulation results for Example 3 in Pin et al. (2017) with different initial values $\tilde{w}(0)$.}
\end{figure*}


\begin{thebibliography}{1}
\bibitem{Angulo}
Angulo M. T., Moreno J. A., \& Fridman L. (2013). Robust exact uniformly convergent arbitrary order differentiator. \emph{Automatica}, 49(8), 2489-2495.

\bibitem{Andrieu}
Andrieu V., Praly L., \& Astolfi A. (2008). Homogeneous approximation, recursive observer and output feedback. \emph{SIAM Journal on Control and Optimization}, 47(4), 1814-1850.

\bibitem{Angrisani}
Angrisani L., D'Apuzzo M., Grillo D., Pasquino N., \& Moriello R. S. L. (2014). A new time-domain method for frequency measurement of sinusoidal signals in critical noise conditions. \emph{Measurement}, 49, 368-381.



\bibitem{Chen1}
Chen B., Pin G., Ng W. M., Hui S. Y. R., \& Parisini T. (2017). An adaptive-observer-based robust estimator of multi-sinusoidal signals. \emph{IEEE Transactions on Automatic Control}, 63(6), 1618-1631.

\bibitem{Chen2}
Chen B., Li P., Pin G., Fedele G., \& Parisini T. (2019). Finite-time estimation of multiple exponentially-damped sinusoidal signals: A kernel-based approach. \emph{Automatica}, 106, 1-7.

\bibitem{Hsu}
Hsu L., Ortega R., \& Damm G. (1999). A globally convergent frequency estimator. \emph{IEEE Transactions on Automatic Control}, 44, 698-713.

\bibitem{Hajimolahoseini}
Hajimolahoseini H., Taban M. R., \& Soltanian-Zadeh, H. (2012). Extended kalman filter frequency tracker for nonstationary harmonic signals. \emph{Measurement}, 45(1), 126-132.

\bibitem{Hou}
Hou M. (2012). Parameter identification of sinusoids. \emph{IEEE Transactions on Automatic Control},  57(2), 467-472, 2012.

\bibitem{Karimi}
Karimi-Ghartemani M., \& Ziarani A. K. (2004). A nonlinear time-frequency analysis method. \emph{IEEE Transactions on Signal Processing}, 52(6), 1585-1595.

\bibitem{Li}
Li P., Fedele G., Pin G., \& Parisini T. (2016). Kernel-based deadbeat parametric estimation of bias-affected damped sinusoidal signals. \emph{In 2016 European Control Conference (ECC)}, 519-524.


\bibitem{Marino}
Marino R., \& Tomei P. (2002). Global estimation of $n$ unknown frequencies. \emph{IEEE Transactions on Automatic Control}, 47(8), 1324¨C1328.


\bibitem{Pyrkin}
Pyrkin A. A., Bobtsov A. A., Efimov D., \& Zolghadri A. (2011). Frequency estimation for periodical signal with noise in finite time. \emph{In 2011 50th IEEE Conference on Decision and Control and European Control Conference}, 3646-3651.

\bibitem{Polyakov}
Polyakov A. (2012). Nonlinear feedback design for fixed-time stabilization of linear control systems. \emph{IEEE Transactions on Automatic Control}, 57(8), 2106-2110.


\bibitem{Pin3}
Pin G., Chen B., \& Parisini T. (2017). Robust finite-time estimation of biased sinusoidal signals: A volterra operators approach. \emph{Automatica},  77, 120-132.

\bibitem{Pin2}
Pin G., Wang Y., Chen B., \& Parisini T. (2019). Identification of multi-sinusoidal signals with direct frequency estimation: An adaptive observer approach. \emph{Automatica}, 99, 338-345.

\bibitem{Rao}
Rao A. V. K., Soni K. M., Sinha S. K., \& Nasiruddin, I. (2019). Accurate phasor and frequency estimation during power system oscillations using least squares. \emph{IET Science, Measurement \& Technology}, 13(7), 989-994.



\bibitem{Shi1}
Shi S., Xu S., Gu J., \& Zhang Z. (2019). Robust exact predictive scheme for output-feedback control of input-delay systems with unmatched sinusoidal disturbances. \emph{IEEE Transactions on Systems, Man, and Cybernetics: Systems}, to be published, doi: 10.1109/TSMC.2019.2951727.

\bibitem{Trapero1}
Trapero J. R., Sira-Ram¨ªrez H., \& Batlle V. F. (2007). An algebraic frequency estimator for a biased and noisy sinusoidal signal. \emph{Signal Process}, 87(6), 1188-1201.


\bibitem{Xia}
Xia X. (2002). Global frequency estimation using adaptive identifiers. \emph{IEEE Transactions on Automatic Control}, 47(7), 1188-1193.


\end{thebibliography}
\end{document}